\input amstex
\documentstyle{amsppt}
\TagsOnRight

\NoRunningHeads\NoPageNumbers\nologo
\footline{}

\catcode`\@=11
\def\raggedcenter@{\leftskip\z@ plus.3\hsize \rightskip\leftskip
 \parfillskip\z@ \parindent\z@ \spaceskip.3333em \xspaceskip.5em
 \pretolerance9999\tolerance9999 \exhyphenpenalty\@M
 \hyphenpenalty\@M \let\\\linebreak}
\catcode`\@=\active

\mag=1200
\pagewidth{14truecm}
\pageheight{21.5truecm}
\hoffset23pt
\voffset36pt
\binoppenalty=10000
\relpenalty=10000
\tolerance=500
\mathsurround=1pt

\define\m1{^{-1}}
\define\ov1{\overline}
\def\gp#1{\langle#1\rangle}
\def\ul2#1{\underline{\underline{#1}}}
\def\supp{\operatorname{supp\,}}
\def\char{\operatorname{char\,}}

\hrule height0pt
\vskip 1truein

\topmatter
\title
On symmetric  units in  group algebras
\endtitle

\author
Victor Bovdi\\
{\eightpoint  University of Debrecen, 
Debrecen,  Hungary\\
}\\
{\eightpoint 
vbovdi\@math.klte.hu
}
\endauthor 

\leftheadtext\nofrills{Victor Bovdi}
\rightheadtext\nofrills{On symmetric  units in  group algebras}

\abstract
Let $U(KG)$ be the group of units of the group ring $KG$ of the
group $G$ over a commutative ring $K$.  The anti-automorphism
$g\mapsto g\m1$ of $G$ can be  extended linearly to an
anti-automorphism $a\mapsto a^*$ of $KG$.  Let $S_*(KG)=\{x\in
U(KG) \mid x^*=x\}$ be the set of all symmetric units of $U(KG)$.
We  consider  the following
question: for which  groups $G$ and commutative rings
$K$ it is  true that $S_*(KG)$  is a subgroup in $U(KG)$.  
We  answer this question when either  
a) $G$ is torsion and $K$ is a commutative $G$-favourable integral 
domain of characteristic $p\geq 0$ or 
b) $G$ is  non-torsion nilpotent group and   
$KG$ is semiprime.
\endabstract


\thanks 
The research was supported by the Hungarian 
National Foundation for Scientific Research Grants No.T 029132  
and No.T 025029.
\endthanks 
\endtopmatter

%\addto\tenpoint{\normalbaselineskip16pt\normalbaselines}

\document

\subhead
1. Introduction
\endsubhead
Let $K$ be a commutative ring, $G$ be a group,  and  $U(KG)$ be
the group of the units of the group ring $KG$.  The
anti-automorphism $g\mapsto g\m1$ of $G$  can be  extended
linearly to an anti-automorphism $a\mapsto a^*$ of $KG$.  An 
element $x\in KG$ is called {\it symmetric}  if $x^*=x$. Let $
S_*(KG)=\{ x\in U(KG)\mid x^*=x\} $ be the set of all symmetric
units of $U(KG)$.  There are  nontrivial symmetric units in
$U(KG)$, for example $xx^*$, where $x\in U(KG)$ and $xx^*\not=1$.

The properties  of the symmetric elements and symmetric units play 
important role in the study of group and ring properties  of rings
 (see  e.g. \cite{2,6,7}).

We say that the commutative ring $K$ is  $G$-favourable if for
any $g\in G$ of finite order there is a non-zero $\alpha\in K$,
such that $1-\alpha^{|g|}$ is a unit in $K$, where $|g|$  denotes
the order of $g$. Every infinite field $K$ is obviously 
$G$-favourable.

The problem when the set of symmetric units form a group was 
first studied  for the case of rings 
in the paper of Lanski \cite{7}, and  for group algebras  in \cite{3,4}.

In the first part of this note we describe those  $K$ and $G$ for which
$S_*(KG)$ is a subgroup of $U(KG)$.  We  obtain the following  
results, which  generalize  a joint result with Kov\'acs and Sehgal
\cite{3}.
 
\proclaim
{Theorem 1}
Let $G$ be a torsion group and  $K$ be a commutative
$G$-favourable integral domain of characteristic $p\geq 0$.  Then
$S_*(KG)$ is a subgroup in $U(KG)$ if and only if  $G$ satisfies
at least one of the following conditions:
\item{\rm{1.}} $G$ is abelian;
\item{\rm{2.}} $p\not=2$ and $G$ is a hamiltonian $2$-group;
\item{\rm{3.}} $p=2$ and $G$ is the direct product of an elementary
abelian group and a group $H$ for which at least one of the
following holds:
\itemitem{\rm{(i)}} $H$ has an abelian subgroup $M$ of index $2$
and an element $b$ of order $4$ such that conjugation by $b$
inverts each element of $M$;
\itemitem{\rm{(ii)}} $H$ is the direct product of a quaternion
group of order $8$ and a cyclic group of order $4$  or the direct
product of two quaternion groups of order~$8$;
\itemitem{\rm{(iii)}} $H$ is the central product of the group
$\gp{\,x,y\mid x^4=y^4=1,\ x^2=[y,x]\,}$  with a quaternion group
of order $8$, the nontrivial element common to the two central
factors being $x^2y^2$;
\itemitem{\rm{(iv)}} $H$ is isomorphic to one of the groups
$H_{32}$ and $H_{245}$ defined below.

\endproclaim

The relevant definitions are:
$$
\align
 H_{32}=\bigl\langle\,x,y,u\bigm|\ &x^4=y^4=1,x^2=[y,x],\\
&y^2=u^2=[u,x], x^2y^2=[u,y]\,\bigr\rangle,\\
\noalign{\vskip2pt}\allowbreak
H_{245}=\bigl\langle\,x,y,u,v\bigm|\ &x^4=y^4=[v,u]=1,
x^2=v^2=[y,x]=[v,y],\\ &y^2=u^2=[u,x],
x^2y^2=[u,y]=[v,x]\,\bigr\rangle.\\
\endalign
$$

In the next part we consider when $S_*(KG)$ is a subgroup for
a  nilpotent groups with element of infinite order and  $KG$ is
semiprime.

Let $t(G)$ be the set of torsion elements of $G$. It is well known  
that for nilpotent group $G$ the set $t(G)$ is a subgroup. 

\proclaim
{Theorem 2}
Let $K$ be an infinite  field of characteristic $p\geq 0$ and $G$
be a nilpotent group such that the characteristic of $K$ does not
divide the order of elements in $t(G)$.  If $S_*(KG)$ is a
subgroup in $U(KG)$ then   $G$  and $K$ satisfy at least one of
the following conditions:
\itemitem{\rm{(i)}} $t(G)$ is a central  subgroup of $G$;
\itemitem{\rm{(ii)}} if  $K$ is  a field of characteristic $0$ then
$t(G)$ is  abelian and for each $x\in G$ there exists a positive
integer $j$ such that $t^x=t^j$ and, if $\phi$ is a primitive
root of unity of order $|t|$, then there exists a map $\sigma\in
Gal(K(\phi):K)$ such that $\sigma(\phi)=\phi^k$;
\itemitem{\rm{(iii)}} If $t(G)$ is  noncentral and $K$ is  a field of
characteristic $p\not=0$, then algebraic closure $\Omega$  of the
prime subfield $P(K)$ in $K$ is finite and for all $x\in G$ and
all $t\in A$ there exists $r\in \Bbb N$, such that
$t^x=t^{p^r}$. Furthermore, for each such $r$  we have
$(\Omega:P(K))|r$;
\itemitem{\rm{(iv)}} $p\not=2$ and $t(G)$ is a hamiltonian $2$-group.
\endproclaim

\subhead
2. Group algebras of torsion groups
\endsubhead
Let $K$ be an integral domain   of characteristic $p\geq 0$. A
group $G$ is called "{\it good}" if $S_*(KG)$ is a subgroup in
$U(KG)$.

To prove our  theorems,  we need the following lemmas. 

\proclaim
{Lemma 1} (cf. \cite{3}) 
If $S_*(KG)$ is a subgroup in $U(KG)$ then $S_*(KG)$ is
abelian and normal in $U(KG)$.
\endproclaim

\rightline{\text{\qed} }

In the sequel  of this section let us assume that $K$ is a 
$G$-favourable integral domain  and $G$
is a torsion group.  If $\alpha\in K$ and
$\alpha^{|g|}-1$ is a unit,  then $g-\alpha\in U(KG)$. 
 Indeed, 
$$
(g-\alpha)\big((\alpha^{|g|-1}+\alpha^{|g|-2}g+\cdots+
{\alpha}g^{|g|-2}+g^{|g|-1})(\alpha^{|g|}-1)\m1\big)=1, 
$$
and it is easy to see  that $(g-\alpha)(g\m1-\alpha)$ is 
a symmetric unit.

Take $g,h,t\in G$ with  $t^2=1$, $g^2\neq1$, $h^2\neq1$, and let
$\alpha$, $\beta$ be the elements corresponding  to $g$, $h$,
respectively, according to  the   definition of
$G$-favourability. Obviously, $t$, $(g-\alpha)(g^{-1}-\alpha)$
and $(h-\beta)(h^{-1}-\beta)$ are  symmetric units, and since  symmetric
units commute, so  $t$,  $\alpha(g+g^{-1})$ and
$\beta(h+h^{-1})$ commute as well. Since $\alpha\beta\neq0$, we
 reached the requived conclusion in  the second  paragraph of
\cite{3} (p.805).  Furthermore, for a $G$-favourable $K$, the
symmetric units of $KG$ form a multiplicative group if and only
if every pair of elements of 
$$ 
S=\{\,\, t\,\, \mid \,\, t\in G ,
\,\,\,t^2=1\,\, \}\cup
\{ \,\, h+h\m1 \,\,\, \mid \,\,\,\,h\in G,\,\,\,  h^2\not=1\,\, \}
$$ 
commute.

\proclaim
{Lemma 2} 
Let $G$ be a good torsion group and let $g,h\in G$, such that 
$[g,h]\not=1$. Then  all involutions of $G$ are central in $G$ 
and:
\itemitem{\rm{(i)}}  if $\char(K)\not=2$,  then either $g^h=g\m1$ or
$h^g=h\m1$ or $\gp{g,h}\cong Q_8$;
\itemitem{\rm{(ii)}}  if $\char(K)=2$,  then $\gp{g, h}=\gp{x, y}$,
where $y^4=1$ and  $x^y=x\m1$. Moreover, the centre
$\zeta(\gp{g,h})$ of $\gp{g,h}$  has exponent $2$ and all of its
elements of odd order form a cyclic subgroup and every
nonabelian $2$-Sylow subgroup is   either a ge\-ne\-ra\-li\-zed
qua\-ter\-ni\-on group or  can be written as 
$$ 
C_{2^m}\rtimes
C_4=\gp{\,\, a, y\,\, \mid \,\, a^{2^m}=y^4=1,\,\,  a^y=a\m1\,\,
}.  
$$ 
\endproclaim 
\demo{ Proof }  
Let $t,h\in G$ with $t^2=1$,  $h^2\not=1$ and  $[t,h]\not=1$.
Then $t, h+h\m1\in S$ so  they  commute. This means   $th+th\m1=ht+h\m1t$,
and from this follows that  $(th)^2=1$.  Hence $th\in S$ and $[th,t]=1$, 
which is impossible.
Therefore, all involutions of $G$ are central.

Let us choose  two elements $g, h\in G$ with $[g,h]\not=1$ and 
$g^2\not=1\not=h^2$. Then  $g+g\m1, h+h\m1\in S$ and we get the 
following equation  
$$
gh+gh\m1+g\m1h+g\m1h\m1=hg+h\m1g+hg\m1+h\m1g\m1.\tag1 
$$ 
We have
to consider the following three cases:

Case 1. If  $gh=h\m1g$ or $gh=hg\m1$, then either $g^h=g\m1$ or  
$h^g=h\m1$.

It is easy to see  that the second case is similar  to the first one, 
so we consider only the first one. 

Set $\char(K)=2$ and $g^h=g\m1$.
Then  $(gh)^2\not=1$ and $gh+(gh)\m1$, $h+h\m1\in S$  commute and 
$$
g(h^2+h^{-2})=g\m1(h^2+h^{-2}).\tag2
$$ 
 
Supposing $g\gp{h}\not=g\m1\gp{h}$, from (2)  we obtain  that
$h^4=1$ and  if  $\gp{g}=\gp{a}\times\gp{b}$, where $|a|=2^t$ and
$|b|$ is odd, then  $a^h=a\m1$ and  $b^h=b\m1$.  As we proved in
\cite{3}, p.806,    $\gp{a, h}$ is  either $C_{2^t}\rtimes C_4$
or a generalized quaternion group.

In the remaining case $g\gp{h}=g\m1\gp{h}$, consequently $\gp{h}$ has
index two   in $\gp{g,h}$, whence  (2) implies $g^2=h^{-4}$.
Therefore,  $gh^2$ is a noncentral involution, which is impossible.

Case 2. If $gh=h\m1g\m1$, then $(gh)^2=1$ and     
$gh$ is central,  whence   $[g, h]=1$, which is impossible.

Case 3. If $gh=g\m1h\m1$, then $g^2=h^{-2}$ and by (1) we have that  
$$  
2gh+gh\m1+g\m1h=2hg+h\m1g+hg\m1.\tag3  
$$  
First, let $\char(K)\not=2$ and $gh=h\m1g$. Then by (3) we have
that 
$gh=h\m1g=hg\m1$ and $g\m1h=gh\m1=hg$,  which imply  $g^2=h^2$,
$h^g=h\m1$, $g^h=g\m1$ and as result  $\gp{g, h}$ is a quaternion
group of order $8$.

If  $\char(K)=2$, then (3) leads to $gh\m1+g\m1h=h\m1g+hg\m1$,
and   then either $gh\m1=hg\m1$ or $gh\m1=g\m1h$.  If
$gh\m1=hg\m1$,   then $(gh\m1)^2=1$ and  $gh\m1$ is  central, which
is impossible.  If $gh\m1=g\m1h$,  then $g^2=h^2=h^{-2} $ and we
can assume that    $x=gh$ and $y=h$. Therefore, $\gp{g,h}=\gp{x, y}$,
where $y^4=1$, $x^y=x\m1$ and   $(xy)\m1=xy\m1\not=xy$.
 Therefore,  Lemma $2$ is  proved.
\enddemo

\rightline{\text{\qed} }

\proclaim
{Corollary}
A good torsion  group $G$ is   locally-finite and 
the set $A$ of elements of odd order of $G$ 
forms an  abelian normal subgroup.
\endproclaim

\demo{ Proof }
Write  $H=\gp{g^2\mid g\in G}$. If $g_1,g_2\in G$, then by Lemma
$2$ we have $[g_1^2,g_2^2]=1$  and $H$ is  normal abelian. Since
$G/H$ is $2$-elementary abelian, $G$ is  locally-finite.
Clearly, the set $A$ is a subgroup in $H$.
\enddemo
\rightline{\text{\qed} }

\proclaim
{Lemma  3}
Let $A$ be a subgroup of odd elements of a  good torsion group
$G$.  If $A$ is not central in $G$,  then  $C_G(A)$ is an abelian
subgroup of  index $2$  and  any $g\in G\setminus C_G(A)$ inverts
all elements in $C_G(A)$.  
\endproclaim

\demo{ Proof }
Let $P$ be a  $2$-Sylow subgroup of $G$ and $P_1=P\cap C_G(A)$.
Suppose  $g\in P\setminus P_1$, $t\in P_1$, $a\in A$ and
$t^2\not=1$, $[g,a]\not=1$. According to  Lemma 2 either 
$g^a=g\m1$ or $a^g=a\m1$.
But  the case  $g^a=g\m1$ is impossible since $a$ is of  odd order.
It is easy to see that $(at)^2\not=1$ and
$[g,at]\not=1$.  Indeed, otherwise $a^{-2}=tgt\m1g\m1\in P$,
which is impossible.  Since  $g+g\m1$ and  $at+(at)\m1$ belong to
$S$, we conclude that 
$$ 
gat+ga\m1t\m1 +g\m1at+g\m1a\m1t\m1=
a\m1t\m1g+a\m1t\m1g\m1+atg\m1+atg.\tag4 
$$ 
If $atg=g\m1a\m1t\m1$,
then $(tg)^2=1$. Thus  $tg\in \zeta(G)$ and $[a,tg]=1$, which is
impossible.  From (4) we have  either $atg=ga\m1t\m1$ or
$agt=g\m1at$.  In the first case, $a^g=a\m1$, it follows that  
$t^g=t\m1$.  If $atg=g\m1at$,  then $a^2\in P$, which is
impossible.

We conclude that $t^g=t\m1$ for all $t\in P_1$ and $P_1$ is
abelian.  Clearly, $C_G(A)=P_1\times A$ is abelian. Since 
every $g\in G\setminus C_G(A)$ inverts all elements of $C_G(A)$, 
we have that  $C_G(A)$ has index $2$ in $G$, and the subgroup   
$P_1$ has index $2$ in $P$.
\enddemo
\rightline{\text{\qed} }

\proclaim
{Lemma  4}
Let $G$ be a nonabelian good torsion group  and let $x,y\in G$ with 
$[x,y]\not=1$.  If $char(K)=p\not=2$ then $\gp{x,y}$ is a 
quaternion subgroup of order $8$, and $G$ is a hamiltonian $2$-group.

\endproclaim
\demo{ Proof }
By Lemma $2$ either   
$\gp{\,\,\,x,y\,\,\,}= \gp{\,\,\, g,h\,\,\,\mid 
     \,\,\,g^{m}=h^{n}=1,  \,\,\, h^g=h\m1}$ or 
$\gp{\,\,\, x,y\,\,\,}= \gp{\,\,\, g,h\,\,\, 
\mid h^{n}=1,g^m\in \gp{h}, \,\,\, h^g=h\m1 }$. 
 Since all involutions of $G$ are central,   $|g|,|h|> 2$.

Let us  consider the first two of the  above  three groups.   
It is easy to see that  $g+g\m1$, $gh+(gh)\m1\in S$ commute and
we have that  
$$ 
(g^2+g^{-2}+2)h^2=g^2+g^{-2}+2.\tag5  
$$ 
It follows that either 
$g^2h^2=g^{-2}h^2$ or  $g^2h^2=1$ or $g^2h^2=g^{-2}$. 
In the first two cases we conclude that $\gp{g,h}\cong Q_8$.

According to the third case  $g^2h^2=g^{-2}$. Then from (5) we get
$$
g^{-6}+2g^{-4}=g^2+2
$$
and $g^4=1$. Thus  $h^2=g^{-4}=1$,  which is a contradiction, so 
$\gp{x,y}$ is a quaternion group of order $8$.  

From  the Corollary, the  elements of odd order form an  abelian subgroup $A$. 
Assume that   $|A|\not=1$.  If $g\in G \setminus
C_G(A)$ is a $2$-element, then for any $a\in A$ 
on the basis of the  first part of the Lemma,  
we obtain that $\gp{a,g}$ is a quaternion group of order $8$,
which is impossible.  Hence, $G=C_G(A)$ and $G=A\times P$, where
$P$ is a Sylow $2$-subgroup of $G$.  Let $\gp{c,d}$ be  a
quaternion subgroup of order $8$ in $P$ and $a\in A$.  Then
$[c,a]=[d,a]=1$ and $[ac,d]\not=1$. As before, we get  that
$\gp{ac,d}$ cannot be  a quaternion subgroup of order $8$. Since $G$ 
is nonabelian,   $|A|=1$ and $G$ is a $2$-group.
 
Let $\gp{c,d}$ be  a quaternion subgroup of order $8$.  Then 
$c^g=c^{\pm 1}$ and
$d^g=d^{\pm 1}$ for any $g\in G\setminus \gp{c,d}$. It   follows
that  $G=Q_8\cdot C_G(Q_8)$.

If $t\in C_G(Q_8)$,  then $\gp{c,dt}$ is a quaternion group of order 
$8$ and  $(td)^c=(td)\m1=\-td\m1$.  
Hence $t^2=1$ and  $C_G(Q_8)$ is an elementary $2$-group.    
Thus  $G$ is  a hamiltonian $2$-group.

\enddemo
\rightline{\text{\qed} }

\demo{ Proof of the 'only if' part of  Theorem 1} 
Let $G$ be a nonabelian torsion group and $K$ be a commutative
$G$-favourable integral domain  of characteristic $p\geq 0$, and
let $S_*(KG)$ be a subgroup.  If $p\not=2$,  Lemma $4$ implies
that $G$ is a hamiltonian $2$-group.

Let  $p=2$  and suppose that $G$ is not  a $2$-group. By Lemma
$3$ the  elements of odd order in $G$ form an abelian subgroup
$A$, and either $A$ is central or $[G:C_G(A)]=2$ and $u^g=u\m1$
for all $g\in G\setminus C_G(g)$ and  $u\in C_G(A)$.  If
$[G:C_G(A)]=2$  then  $M=C_G(A)$ is  abelian  and there exists 
 $b\in G$  such that the conjugation by $b$ inverts each element 
of $M$. Since all involutions are central, $b$ has order $4$  and $G$
satisfies  condition \rm{3.(i)} of the   Theorem.

Assume now that  $A$ is   central in  $G$. Then
$G=A\times P$ and according to  \cite{3} the $2$-Sylow subgroup $P$ is a
direct product of an elementary $2$-group and one of those  $H$ 
which were mentioned  in the conditions \rm{3.(i)--3.(iv)}  of  
Theorem 1.

We will prove that $|A|=1$. Indeed, it is easy to see that $P$
contains the subgroup $\gp{a,b}$ such that  $a^4=b^4=1$ and
$a^b=a\m1$.  If $t\in A$ is an element of odd order and
$\gp{a,b}\subset P$ is not the quaternion group of order $8$,
then the  symmetric units 
$x=1+(a+b)t+(a\m1+b\m1)t\m1$ and  
$y=1+(a+b\m1)t+(a\m1+b)t^{-1}$ 
 noncommute in $U(F\gp{a,b})\subset U(KG)$, where $F=KA$.  
If $t\in A$
and $\gp{a,b}\subset P$ is the  quaternion group of order $8$,
then the symmetric units  
$x=1+(t\m1+ta^2)(a+b)$ and 
$
y=1+(t\m1+ta^2)a(1+b) 
$ 
noncommute either, which is impossible. Hence $|A|=1$ and $G$
satisfies one of the conditions \rm{ 3.(i)--3.(iv)}  of  Theorem 1.
\enddemo
\demo{ Proof of the 'if' part of  Theorem 1} 
Let $p\not=2$ and $G$ be a finite hamiltonian $2$-group. Then
$G=Q\times E$, where $Q=\gp{a,b\mid a^4=1, a^2=b^2, a^b=a\m1}$ 
and $E$ is  $2$-elementary abelian. If 
$H=\gp{a^2}\times E$, then  any $x\in U(KG)$ can be written as
$$
x=x_0+x_1a+x_2b+x_3ab,
$$ 
where $x_i\in KH$ and $x^*=
x_0^*+(x_1a+x_2b+x_3ab)a^2$. Clearly $x\in S_*(KG)$ if and only
if $x_0$ is central and $(x_1a+x_2b+x_3ab)(1+a^2)=0$. We know
that $x_1a+x_2b+x_3ab=(1+a^2)y$ with  some $y\in KG$ and all
elements of this form are central. Therefore,  all symmetric units
are central in $KG$ and $S_*(KG)$ is a subgroup.

Suppose  that $p=2$ and $G$ has an abelian subgroup $M$ of index
$2$ and an element $b$ of order $4$ such that  conjugation by
$b$ inverts each element of $M$.  Any  $x\in U(KG)$  can be
written as $x=x_0+x_1b$ and $x^*=x_0^*+x_1b^2b$, where
$x_0,x_1\in KM$.  Clearly, $x\in S_*(KG)$ if and only if $x_0$ is
 central  and $x_1(1+b^2)=0$. It  follows that $x_1=(1+b^2)x_2$. 
Let $x=x_0+x_1b, y=y_0+y_1b\in S_*(KG)$, where 
$x_1=(1+b^2)x_2, y_1=(1+a^2)y_2$. Since
$$
x_1y_1^*=x_2(1+a^2)y_2^*(1+a^2)=0
$$ 
and similarly $x_1^*y_1=0$, we
obtain that 
$$ 
xy=(x_0+x_1b)(y_0+y_1b)=x_0y_0+(x_1y_0+x_0y_1)b=yx
$$ 
and $S_*(KG)$ is a subgroup.

If $G$ satisfies  one of conditions \rm{3.(ii)--3.(iv)} of the Theorem,
then by (3) $S_*(KG)$ is a subgroup of $U(KG)$.

\enddemo
\rightline{\text{\qed} }

\subhead
3. Group algebra of  a  nilpotent groups with element of infinite order 
\endsubhead

\proclaim
{Lemma  5}
(see Lemma 1 in \cite{1})
Suppose that $G$ is a nilpotent  group and the characteristic of
$K$ does not divide the order of elements in $t(G)$.  Let
$x_1,\ldots,x_s $ be finite number  elements in $U(KG)$ and
assume that all idempotents of $Kt(G)$ are central in $KG$.  Then
there exists a finite subgroup $L$ of $t(G)$ such that 
$$
x_j=\sum_{i=1}^t\alpha_ie_ig_i\,\,\,
\text{ and }\,\,\,
x_j\m1=\sum_{i=1}^t\alpha_i\m1e_ig_i\m1,
$$
where the list $e_1,\ldots,e_t$ contains  all primitive
orthogonal idempotents  of $KL$, $\alpha_i\in U(KLe_i)$ and
$g_i\in G$ for all $j=1,\ldots,s$.
\endproclaim

\rightline{\text{\qed} }

\subhead
4. Proof of  Theorem 2
\endsubhead
Suppose that $S_*(KG)$ is a subgroup in $U(KG)$ and the
characteristic of $K$ does not divide the order of elements in
$t(G)$, where $t(G)$ is the set of torsion elements of $G$. 
Then  $S_*(Kt(G))$ is a subgroup of $U(Kt(G))$ and by
Theorem 1, either $t(G)$ is abelian or  $p\not=2$  and $t(G)$ is
a  hamiltonian $2$-group.

Let us  assume that $t(G)$ is an abelian group. By \cite{5} all
idempotents of $Kt(G)$ are central in $KG$ if and only if the
conditions (i)--(iii) of  Theorem 2 are satisfied. 

Evidently,  to prove Theorem 2   it is enough to show   
that all idempotents of $Kt(G)$ are central in $KG$.

Now, let $f^2=f\in Kt(G)$ and suppose that there exists 
$g\in G$ such that $[f,g]\not=1$. It is easy to see that
$L=t(\gp{g,\supp(f)})$ is a finite abelian  group. According to the  
Wedderburn-Artin theorem $e_1+e_2+\cdots+e_s=1$ and
$f=e_{i_1}+\cdots +e_{i_k}$, where $e_1,\ldots, e_s$ are
primitive orthogonal idempotents of $KL$. It follows that there
exists a primitive idempotent $e\in KL$ such that $[e,g]\not=1$. 
If $e\not=e^*$, then $e^*=e_j$ for some $j=1,\ldots,t$. Indeed,
suppose that  $e^*=e_{i_1}+\cdots +e_{i_m}$. Then $e=e_{i_1}^*+\cdots
+e_{i_m}^*$. Since $e_{i_k}e_{i_l}=0$  $(k\not=l)$, we have
$e_{i_k}^*e_{i_l}^*=0$ $(k\not=l)$ and $e_k$ is a sum of
orthogonal idempotents, which contradicts the fact that $e_k$ 
is  primitive.

We have  six  cases to consider:

Case 1. Let $e=e^*$. Then 
$e_{11}=e$,
$e_{12}=eg$, 
$e_{21}=g\m1e$ 
and 
$e_{22}=g\m1eg$ 
are matrix units. If $f=e_{11}+e_{22}$ and 
$
A=\left(\smallmatrix
\alpha_1 & \alpha_2 \\ 
 \alpha_3 & \alpha_4
\endsmallmatrix\right)
\in GL(2,K)$, then
$$
w(A)=1-f+\alpha_1e_{11}+\alpha_2e_{12}+\alpha_3e_{21}+ 
\alpha_4e_{22}\in U(KL).
$$
Indeed, $w(A)w(A\m1)=1-f+e_{11}+e_{22}=1$.

Clearly $w(A)\in S_*(KG)$ if and only if $\alpha_2=\alpha_3$. 

Let 
$
A=\left(\smallmatrix
1 & t\\ t & t
\endsmallmatrix\right)
$
and 
          $ 
          B=\left(\smallmatrix
          t   & t+1 \\ 
          t+1 & 1
          \endsmallmatrix\right)
          $
are  units in $GL(2,K)$, such that $t$ does  not satisfy
the equation $2x^2-x-1=0$ in $K$. Since $K$ is infinite,  such
$x\in K$ always exists.  We have  $AB\not=BA$ and  $w(A),w(B)\in
S_*(KL)$.  Therefore, $w(A)w(B)\not=w(B)w(A)$, which is impossible
by Lemma $1$.

Case 2. Let  $e=e_1$ and let $e_1^*=e_2$, $e_1^g=e_2$, $e_1^{g^2}=e_1$. 
Then  
$e_{11}=e_1$, $e_{12}=e_1g$, $e_{21}=g\m1e_1$ and
$e_{22}=g\m1e_1g=e_2$  are matrix units and 
$
e_{11}^*=e_{22}$, $e_{12}^*=g^{-2}e_1g=g^{-2}e_{12}$,
$e_{21}^*=g\m1e_1g^2=e_{21}g^2$, $e_{22}^*=e_{11}$.
  
Let 
$ A=\left(\smallmatrix
\alpha_1 & \alpha_2 \\ 
 \alpha_3 & \alpha_4
\endsmallmatrix\right)
\in GL(2,K\gp{g^2})$ 
and $f=e_{11}+e_{22}$. Then
$$
w(A)=1-f+\alpha_1e_{11}+\alpha_2e_{12}+ \alpha_3e_{21}+
\alpha_4e_{22}\in U(KL).
$$ 
Since $[g^2,e_1]=[g^2,e_1^*]=1$, we have $w(A)\in S_*(KL)$ if and
only if $\alpha_1=\alpha_4^*$, $\alpha_2=\alpha_2g^{-2}$,
$\alpha_3=\alpha_3g^{2}$, $\alpha_4=\alpha_1^*$.

If  
$
A=\left(\smallmatrix
\format\l\quad&\l\\
g^2 & g^{-2}+1 \\ 
0   & g^{-2}
\endsmallmatrix\right)
$ and
$
B=\pmatrix
\format\l\quad&\l\\
g^6      & 0 \\
g^{2}+1 & g^{-6}
\endpmatrix
$, 
then $w(A)$ and  $w(B)$ belong to $S_*(KL)$, but 
$$
AB=
\pmatrix
\format\l\quad&\l\\
2+g^8+g^2+g^{-2} & g^{-8}+g^{-6} \\ 
1+g^{-2}   & g^{-8}
\endpmatrix
\not=
$$$$
\not=BA=
\pmatrix
\format\l\quad&\l\\
g^8      &  g^4+g^{6} \\ 
1+g^4    &  2+g^{-8}+g^2+g^{-2}
\endpmatrix
$$ 
and $w(A)w(B)\not=w(B)w(A)$, which is impossible
by Lemma $1$.

Case 3. Let $e=e_1$, $e_1^*=e_2$, $e_1^g=e_2\not=e_2^*$ and
$e_2^g=e_3\not=e_1$. Then
$$
e_3^*=(g\m1e_2g)^*=g\m1e_2^*g=g\m1e_1g=e_2,
$$
whence  $e_3=e_2^*=e_1$, which is impossible.

Case 4. Let $e=e_1$, $e_1^*=e_2$, $e_1^g=e_3\not=e_2$ and
$e_2^g=e_2$. Then
$$
e_3^*=(g\m1e_1g)^*=g\m1e_1^*g=g\m1e_2g=e_2=e_1^*
$$
and $e_1=e_3$, which is impossible.

Case 5. Let $e=e_1$, $e_1^*=e_2$, $e_1^g=e_3\not=e_2$ and
$e_2^g=e_3$. Then $e_3^*=g\m1e_2^*g=g\m1e_1g=e_3$  and we have
Case 1.

Case 6. Let $e=e_1$, $e_1^*=e_2$, $e_1^g=e_3\not=e_2$ and
$e_2^g=e_4\not=e_3$. Denote  $\ov1{e}=e_1+e_2$. Then 
$\ov1{e_1}^*=\ov1{e_1}$ and  we have Case 1.

Therefore, all idempotents of $Kt(G)$ are central in $KG$ and  this 
completes the proof of Theorem $2$.

\rightline{\text{\qed} }

The question of when these conditions are sufficient
to prove that  symmetric units form a subgroup  is still open.

Let us  make the following remarks for the case when  $G$ is
nilpotent and $t(G)$ is abelian:

{\bf Remark 1.} Suppose that all idempotents $e$ of $Kt(G)$ are
symmetric, i.e. $e^*=e$.  By Lemma $5$  any $x\in U(KG)$ can be
written as $x=\sum_{j=1}^t\alpha_je_jg_j$.  Then the unit $x$ is
symmetric if and only if 
$$
xe_k=\alpha_kg_ke_k=g_k\m1\alpha_k^*e_k=xe_k 
$$ 
for all $k$. It
follows that $g_k$ is an  element of finite order.  Therefore, if
$x\in S_*(KG)$,  then $x\in U(Kt(G))$. Since $ U(Kt(G))$ is
abelian, we obtain that $S_*(KG)$ is a subgroup.

{\bf Remark 2.} Suppose that $Kt(G)$ contains a non-symmetric
idempotent $e$.  We show that in  this case  there exists a
non-symmetric idempotent $f$ such that $f$ and $f^*$ are
orthogonal.

Let $H$ be a  supplement subgroup of $e$. Then $H$ is a finite
abelian subgroup and the group algebra $KG$ is semisimple.
Clearly, every idempotent of $KH$ is a sum of primitive
idempotents and let $e=e_1+e_2+\ldots+e_s$, where $e_i$ is a
primitive idempotent. Because $e\ne e^*$, there exists  $e_i$
which is not symmetric. We show that if $e_i$ is a  primitive
idempotent, then $e_i^*$ is also  a primitive idempotent.  Indeed,
let $e_i^*=e_{i_1}+\cdots+e_{i_p}$ be  a sum of primitive
idempotents. Then $e_i=e_{i_1}^*+\cdots+e_{i_p}^*$. Since
$e_{i_k}e_{i_l}=0$ $(k\not=l)$,  we have $e_{i_k}^*e_{i_l}^*=0$
(${k\not=l}$) and $e_i$ is a sum of orthogonal idempotents, which
contradicts the fact that $e_i$ is  primitive.

{\bf Remark 3.} Suppose that $K$ is  a field with the following
property: the equation $x_1^2+x_2^2+\cdots+x_s^2=0$ has only
trivial solution  in $K$. Then all idempotents are symmetric.
Indeed,  there exists a non-symmetric idempotent $f=\sum_{g\in G}
\alpha_g g$  such that $f$ and $f^*$ are orthogonal and
$0=tr(ff^*)=$ $=\sum_{g\in G}\alpha^2$, which is impossible.

{\bf Remark 4.} Suppose that all idempotents of $Kt(G)$ are
central in $KG$ and $Kt(G)$ contains a non-symmetric idempotent.
Choose a non-symmetric idempotent $f$ such that $f$ and $f^*$ are
orthogonal. Then every  $x_g=1-f-f^*+gf+g^{-1}f^*$ is a   symmetric
unit for every $g\in G$. Suppose that all
$x_gx_h=1-f-f^*+ghf+g^{-1}h^{-1}f^*$ are  symmetric units. Then
$(gh)^{-1}f^*=g^{-1}h^{-1}f^*$. It implies that $((g,h)-1)f^*=0$
and the commutator subgroup of $G$ belongs to the supplement of
the idempotent $f^*$. Therefore, if $KG$ has non-symmetric units,
the commutator subgroup is a finite group.

{\bf Remark 5.}  
Let $char(K)\not| 2$ and let $G=Q\times W$  be a hamiltonian $2$-group,   
where $Q=\gp{a,b}$ is a quaternion group  and $exp(W)=2$.

Clearly, $KG=RQ$, where $R=KW$ and  $RQ=Re_1\oplus\cdots\oplus Re_5$, 
where 
\newline
$e_1=\frac{1}{8}(1+a)(1+a^2)(1+b)$,
$e_2=\frac{1}{8}(1-a)(1+a^2)(1+b)$,
\newline
$e_3=\frac{1}{8}(1+a)(1+a^2)(1-b)$,
$e_4=\frac{1}{8}(1-a)(1+a^2)(1-b)$,
$e_5=\frac{1}{2}(1-a^2)$,
and $Re_i\cong  R$, for $i=1,\ldots,4$.

It is easy to see that any $x\in Re_i$  ($i=1,2,3,4$) is 
symmetric element. The element $x\in Re_5$ can be written as
$x=(\alpha_0+\alpha_1a+\alpha_2b+\alpha_3ab)e_5$ and $x=x^*$ 
if and only if $\alpha_1=\alpha_2=\alpha_3=0$. Therefore, 
all symmetric elements of $KG$ are central and $S_*(KG)$ 
is a subgroup.

\enddocument